\documentclass[11pt]{article}

\usepackage{geometry}
\usepackage{amsmath}
\usepackage{amssymb}
\usepackage{amsthm}
\usepackage[latin1]{inputenc}
\usepackage{eurosym}
\usepackage[dvips]{graphics}
\usepackage{graphicx}
\usepackage{epsfig}
\usepackage{hyperref}

\usepackage{ifthen}

\makeindex

\newcommand{\Rr}{{\mathbb{R}}}

\newcommand{\Tt}{{\mathbb{T}}}

\newcommand{\Gg}{{\mathcal{G}}}

\theoremstyle{definition}

\newtheorem{theorem}{Theorem} \newtheorem{corollary}{Corollary}

 \newcommand{\re}{\mathbb{R}}
\newcommand{\tn}{\mathbb{T}^N}
\newcommand{\rn}{\mathbb{R}^N}

\begin{document}

\title{Exponential decay of correlation for the Stochastic Process associated
to the Entropy Penalized Method }
\author{D. A.  Gomes (IST - Portugal) and A. O. Lopes (Inst. Mat - UFRGS -Brasil)}
\date{\today}

\maketitle

\thanks{D. Gomes was
partially supported by the Center for Mathematical Analysis,
Geometry and Dynamical Systems through FCT Program POCTI/FEDER and
also by grant POCI/FEDER/MAT/55745/2004. A. O. Lopes \, was \,
partially \, supported by CNPq, PRONEX -- Sistemas Din\^amicos,
Instituto do Mil\^enio, and is beneficiary of CAPES financial
support. }

\begin{abstract}
In this paper we present an upper bound for the decay of correlation  for the
stationary stochastic process associated with
the Entropy Penalized
Method.

Let $L(x, v):\Tt^n\times\Rr^n\to \Rr$ be a Lagrangian of the form
\[
L(x,v) \,=\, \frac{1}{2}\, |v|^2\,\, - \,U(x) \,+\, \langle P, v\rangle.
\]
For each value of $\epsilon $ and $h$, consider the
operator
$$\Gg[\phi](x):= -\epsilon h \mbox{ ln}\left[
  \int_{\re^N} e ^{-\frac{hL(x,v)+\phi(x+hv)}{\epsilon h}}dv\right] ,$$
as well as the reversed operator
$$\bar \Gg[\phi](x):= -\epsilon h \mbox{ ln}\left[
  \int_{\re^N} e ^{-\frac{hL(x+hv,-v)+\phi(x+hv)}{\epsilon
      h}}dv\right],$$
both acting on continuous functions $\phi:\Tt^n\to \Rr$.
Denote by $\phi_{\epsilon,h} $ the solution of
$\Gg[\phi_{\epsilon,h}]=\phi_{\epsilon,h}+\lambda_{\epsilon,h}$,
and
by $\bar \phi_{\epsilon,h} $ the solution of
$\bar \Gg[\phi_{\epsilon,h}]=\bar \phi_{\epsilon,h}+\lambda_{\epsilon,h}$.
Let
$ \theta_{\epsilon,h}(x)=
 e^{-\, \frac{\bar  \phi_{\epsilon,h}(x)+ \phi_{\epsilon,h}(x)   }{\epsilon \, h }} $.
From \cite{GV}, it is known that
$$\mu_{\epsilon,h}(x,v)=\theta_{\epsilon,h}(x) \,\gamma_{\epsilon,h}(x,v)\,= \,
\theta_{\epsilon,h}(x) \,  e^{-\frac{hL(x,v)+\phi_{\epsilon,h}(x+hv)-\phi_{\epsilon,h}(x)-\lambda_{\epsilon,h}}{\epsilon h}},$$
is a solution to the entropy penalized problem:
$\min \{ \int_{\tn\times\rn} L(x,v)d\mu(x,v)+\epsilon S[\mu]\},$
where the entropy $S$ is given by
$$S[\mu]=\int_{\tn\times\rn} \mu(x,v)\ln\frac{\mu(x,v)}{\int_{\rn}\mu(x,w)dw}dxdv,
$$
and the minimization is made over all holonomic probability
densities on $\Tt^n\times \Rr^n$, that is probabilities that
satisfy $\int \varphi(x+v)-\varphi(x) \mu(x,v) dxdv=0$, for all
$\varphi \in C^1 (\Tt^n)$. The density $\gamma_{\epsilon,h}(x,v)$
defines a Markovian transition kernel on
$(\mathbb{T}^n)^\mathbb{N}$. The invariant initial density in
$\mathbb{T}^n $ is $\theta_{\epsilon,h}(x).$ In order to analyze
the decay of correlation for this process we show that the
operator $  \, {\cal L} \, (\varphi)\, (x) = \int \, e^{-\, \frac{
h L ( x,v) }{\epsilon  }   }   \, \varphi(x+h\, v) \, d \, v,$ has
a maximal eigenvalue isolated from the rest of the spectrum.
\end{abstract}

\section{Definitions and the set up of the problem}

Let $\mathbb{T}^n$ be the $n$-dimensional
 torus.
In this paper
we assume that the Lagrangian, $L(x,v):\Tt^n\times \rn\to \Rr$
 has the form
$$L(x,v) \,=\, \frac{1}{2}\, |v|^2\,\, - \,U(x) \,+\, \langle P, v\rangle,$$
where
$U\in C^\infty(\Tt^n)$,
and $P\in \Rr^n$ is constant.

We consider here the discrete time Aubry-Mather problem \cite{Gom}
and the Entropy Penalized Mather method which provides a way to
obtain approximations by continuous densities of the Aubry-Mather
measure. We refer the reader to \cite{Gom} and the last section of
\cite{GLM} for some of  the main properties of Aubry-Mather
measures , subactions, Peierl's barrier, etc...

The Entropy Penalized Mather problem (see \cite{GV} for general
properties of this problem) can be used to approximate
 Mather measures \cite{CI}  by means of absolutely continuous densities
$\mu_{\epsilon,h}(x)$, when $\epsilon,h\, \to 0$, both in the continuous
case or in the discrete case. In \cite{GLM} it is presented a
Large Deviation principle associated to this procedure.
%
%
We briefly mention some definitions and results.

\noindent {\bf Definition 1:} The forward (non-normalized)
Perron operator ${\cal L}$ is defined
$$ \,\, x \, \to \varphi(x) \, \Rightarrow \, x \to  \, {\cal L} \, (\varphi)\, (x) = \int \,
e^{-\, \frac{  L ( x,v)    }{\epsilon  }   }   \, \varphi(x+h\, v)
\, d \, v,$$

\medskip

In \cite{GV} it is shown that ${\cal L}$ has a unique eigenfunction
$e^{-\frac{\phi_{\epsilon,h}}{h \,\epsilon}  } $
with eigenvalue  $e^{-\frac{\lambda_{\epsilon,h}}{h \,\epsilon}  } $

\medskip

\noindent {\bf Definition 2:} The backward operator ${\cal N}$ is given by
$$ \,\, x \, \to \varphi(x) \, \Rightarrow \, x \to  \, {\cal N} \, (\varphi)\, (x) = \int \,
e^{-\, \frac{  L ( x- h v,v)    }{\epsilon  }   }   \,
\varphi(x-h\, v) \, d \, v,$$

\medskip

In \cite{GV} it is shown that ${\cal N}$ has a unique eigenfunction
$e^{-\frac{\bar \phi_{\epsilon,h}}{h \,\epsilon}  } $
with eigenvalue  $e^{-\frac{\lambda_{\epsilon,h}}{h \,\epsilon}  } $

\medskip

\noindent {\bf Definition 3:}
The operator
$$g(x) \to {\cal F} (g)\, (x)  =  \int
 e^{-\frac{hL(x,v)+\phi_{\epsilon,h}(x+hv)-\phi_{\epsilon,h}(x)-\lambda_{\epsilon,h}}{\epsilon h}}
\, g(x+h v) \, dv,$$
is the normalized forward  Perron  operator.

\medskip

From \cite{GV} we have that given a continuous function
$g:\Tt^n\to \Rr$, then ${\cal F}^m (g)$ converges to the unique
eigenfunction $k$ as $m\to \infty$.
We show in this paper that for $\epsilon$ and $h$ fixed,
the convergence is exponentially fast.

Denote by
$\theta=  \theta_{\epsilon,h}(x)= e^{-\, \frac{\bar  \phi_{\epsilon,h}(x)+ \phi_{\epsilon,h}(x)   }{\epsilon \, h }} $,
$$\gamma(x,v)=\gamma_{\epsilon, h}(x,v)\,
= e^{-\frac{hL(x,v)+\phi_{\epsilon,h}(x+hv)-\phi_{\epsilon,h}(x)-
\lambda_{\epsilon,h}}{\epsilon h}},$$
and set $\mu_{\epsilon,h}=\theta_{\epsilon,h}(x)\gamma_{\epsilon,h}(x, v)$.
Then $\theta_{\epsilon,h}(x)$ is a probability measure in $\Tt^n$ (up to the
addition of a suitable constant to $\phi_{\epsilon,h}(x)$ and
$\bar \phi_{\epsilon,h}(x)$), $\gamma_{\epsilon,h}(x,v)$ is, for each $x\in \Tt^n$,
a
probability measure in $\Rr^n$ and $\mu_{\epsilon,h}(x,v)$ is a probability
measure on $\Tt^n\times \Rr^n$.
A measure $\nu(x,v)$ on $\Tt^n\times \Rr^n$ is called holonomic if
$$\int_{ \Tt^n\times \Rr^n    }\left[\frac{\varphi(x+hv)-\varphi(x)}{h}\right]d\nu(x,v)=0,
$$ for any continuous $\varphi(x)$.
The measure $\mu_{\epsilon, h}$ is an example of an holonomic measure.

\section{Reversed Markov Process and Adjoint Operator}

In this section we define the reversed Markov process and compute
the adjoint of ${\cal F}$  in ${\cal L}^2 (\theta)$. We assume
$h=1$   from  now on.

We can consider the stationary forward Markovian process  $X_n$ according to
the initial probability  $\theta(x)$ and transition $\gamma(x,v)$.
For example
\[
P(X_0\in A_0)=\int_{x \in\mathbb{T}^n\cap A_0}  \theta(x) dx,
\]
$$ P(X_0\in A_0, X_1 \in A_1) \,=\, \int_{x \in\mathbb{T}^n\cap A_0,\,\,( x+ v) \in A_1 }
\theta(x) \gamma(x,v) \, dx \, dv,$$
and so on.
Define the backward transfer operator ${\cal F}^* $ acting
on continuous functions $f(x)$
by
$${\cal F}^* \,  (f)\, (x) = \int \frac{\theta (x-v)   \, \gamma(x-v,v)}{\theta(x)}\, f(x-v) \, dv.                          $$
The backward transition kernel is given by
$$Q(x,v)=\frac{\theta (x-v)   \, \gamma(x-v,v)}{\theta(x)}.$$
The fact that for any $x$ we have  $\int Q(x,v)\,\, dv =1$ follows
from Theorem 32 in \cite{GV}.
We will show in Corollary \ref{cor1}
that $\theta$ is an invariant measure
for the process with transition kernel $Q$, more precisely,
that
\[
\int g d\theta=\int {\cal F}^*(g) d\theta,
\]
for any $g\in {\cal L}^2(d\theta)$.

\begin{theorem}
${\cal F}^* $ is the adjoint of ${\cal F}$  in ${\cal L}^2 (\theta)$,
that is for all $f, g\in {\cal L}^2(\theta)$ then
\[
\int f(x) {\cal F} g(x)\theta(x) dx =
\int g(x) {\cal F}^* f(x)\theta(x) dx.
\]
\end{theorem}
\begin{proof}
Consider $f,g \in {\cal L}^2 (\theta)$, then

\begin{align*}
\int\,  g(x)\,\,[\,  {\cal F}^*  (f)\, (x) \,]\,
\theta (x) dx \, &= \,  \int\,  g(x)\,\,
[\,\int \frac{\theta (x-v)   \, \gamma(x-v,v)}{\theta(x)}\, f(x-v) \,
dv\,] \,\,\theta (x) \,\, dx  \\
&=
\int\,  g(x)\,\, [\,\int \, \theta (x-v)   \, \gamma(x-v,v) \, \, f(x-v) \, dv\,]  \,\, dx  \\& =
\int\, [\, \int \, [ g(x)\,\, \theta (x-v)   \, \gamma(x-v,v) \, \, f(x-v) \,] \, dx\, ] \,\, dv \\&  =
\int\, [\, \int \,  g(x+ v)\,\,  \,\theta (x)   \, \gamma(x,v) \, \, f(x) \, dx\,]  \,\, dv  \\& =
\int\,  f(x)  \, [\,\int    \, \gamma(x,v) \,\, g(x+ v) \, dv\,]  \,\, \theta(x)\,  dx   \\ &=
\int\,  f(x)\,\, [\,\int
 e^{-\frac{L(x,v)+\phi_{\epsilon,1}(x+v)-\phi_{\epsilon,1}(x)-\lambda_{\epsilon,1}}{\epsilon }}
\, g(x+ v) \, dv\,]  \,\, \theta(x)\,dx \,\\&=
\, \int f(x) \, [\, {\cal F}(g) \,(x)\,] \, \theta(x)\,  dx,
\end{align*}
where we use above the change of coordinates $x \to x - v$ and the fact
that $\mu$ is holonomic.
\end{proof}

\begin{corollary}
\label{cor1}
Consider
the inner product $\langle\cdot,\cdot\rangle$
in ${\cal L}^2 (\theta)$.
Then ${\cal F}$ leaves invariant the orthogonal space to the constant
functions:
$\{g \,|\,\langle g,1\rangle\,=\,   \int\, g\, 1\, d \theta=0\}$.
Furthermore
\[
\int g d\theta=\int {\cal F}^*(g) d\theta.
\]
\end{corollary}
\begin{proof}
Note that ${\cal F} (1)=1$, therefore
\[
\int g 1 d\theta= \int g {\cal F}(1) d\theta= \int {\cal F}^*(g) d\theta.
\]
Thus if $\int g 1 d\theta=0$ it follows $\int {\cal F}^*(g) d\theta=0$.
%
%
%
\end{proof}

\section{Spectral gap, exponential convergence and decay of correlations}

From \cite{GV} it is known that
${\cal L}$ has a unique (normalized) eigenfunction
$e^{-\frac{\phi_{\epsilon,h}}{h \,\epsilon}  } $ corresponding to the largest
eigenvalue  $e^{-\frac{\lambda_{\epsilon,h}}{h \,\epsilon}  } $, in the next
theorem we prove the this eigenvalue is separated from the rest of the
spectrum.

\begin{theorem}
The largest  eigenvalue of ${\cal L}$
is at a positive distance from the rest of the spectrum.
\end{theorem}
\begin{proof}
We will prove the result for the normalized operator
$$  g(x) \to \,{\cal F} (g)\, (x)  =  \int
 e^{-\frac{hL(x,v)+\phi_{\epsilon,h}(x+hv)-\phi_{\epsilon,h}(x)-\lambda_{\epsilon,h}}{\epsilon h}}
\, g(x+\,h\, v) \, dv.$$

Recall from \cite{GV} that the functions $\phi_{\epsilon,h}(x)$
and $\bar \phi_{\epsilon,h}(x)$ are differentiable. In this way we
consider a new Lagrangian (adding
$\phi_{\epsilon,h}(x+hv)-\phi_{\epsilon,h}(x)-\lambda_{\epsilon,h}$)
in such way ${\cal L}={\cal F}$. We also assume $\epsilon=1$ and
$h=1$ from now on.

Therefore,
$$  g(x) \to {\cal F} (g)\, (x)  =  \int e^{- L(x,v)}
\, g(x+v) \, dv,$$ the eigenvalue is $1$, and, by the results in
\cite{GV}, the corresponding eigenspace is one-dimensional and is
generated by the constant functions.

Suppose there exist a sequence of $f_p \in {\cal L}^2 (\theta )$,
$p \in \mathbb{N}$. such that
$${\cal F} (f_p) = \lambda_p (f_p),$$ $\langle f_p, 1\rangle=0$,
$\lambda_p \to 1$ and $||f_p|| =1$. If the operator is compact,
then the theorem follows from the classical argument: through a
subsequence $f_p\to f$, and since $\lambda_p\to 1$ we have ${\cal
F} (f)=f$. Furthermore, since $\langle f_p,1\rangle=0$, it follows
$\langle f, 1\rangle=0$, which is a contradiction. Therefore we
proceed to establish the compactness of the operator $\cal F$.


To establish compactness, consider
 $g \in {\cal L}^2 (\theta )  $. We claim that
$f={\cal F} (g)$ is in the Sobolev space
${\cal H}^1$ (see \cite{E} for definition and properties).
Indeed, for a fixed $x$, we will compute the derivative of $f$.
Integrating by parts we have
\begin{align*}
\frac{d}{dx} \, f(x) &=  \frac{d}{dx}\, ( {\cal F} (g)\, (x)) \\&=
\int    (\,
 [ \frac{d}{dx} g(x+v) ]\,\,    e^{- L(x,v)}
 \,  - L(x,v) \,[\, \frac{d}{dx}
 \,     e^{- L(x,v)} \,  ]\,\, g(x+v)
\, ) \, dv \\
& = \int (\,
 [\, \frac{d}{dv} g(x+v) \,]\,\,    e^{- L(x,v)}
 \,   -L(x,v) \,[\,\frac{d}{dx}
    e^{- L(x,v)}   \,]\,\, g(x+v)
\,) \, dv \\ &= \int
(\, [ \frac{d}{dv}  e^{- L(x,v)} \,]\, g(x+v) \,\,
  -L(x,v) \,[\,\frac{d}{dx}
\, e^{- L(x,v)}   \,]\,\, g(x+v) \,) \, dv  \\& = \int
(\, [ \frac{d}{dv}  e^{- L(x,v)} \,]\,  \,\,
  -L(x,v) \,[\,\frac{d}{dx}
\, e^{- L(x,v)}
\,]\,)\, g(x+v) \, \, dv  .
\end{align*}
From the hypothesis about $L$, if
$g\in {\cal L}^2(\theta)$, then indeed $\frac{d}{d x}\,f$ is also in
${\cal L}^2 (\theta)$
(with the above derivative).

Note that, for $v$ uniformly in a bounded set
$$ \left\| \frac{d}{dx} \, f
\right\|_2\leq
\left\| \frac{d}{dx}\, f\right\|_\infty \leq \left\|\,
[\, \frac{d}{dv}  e^{- L(x,v)} \,]\,  \,\,
  -L(x,v) \,[\,\frac{d}{dx}
\, e^{- L(x,v)}
\,]\,\right\|_2 \, \| g\|_2.$$
Therefore, $f$ is in the Sobolev space ${\cal H}^1$.

By iterating the procedure described above, we have that
$$g_j={\cal F}^j (g) \in {\cal H}^j.$$

It is known that if $j>\frac{n}{2}$, where $n$ is the dimension of
the torus $\mathbb{T}^n$, then $g_j$ is continuous H\"older
continuous\cite{E}. Thus the operator ${\cal F}$ is compact and
$g_j$ is differentiable for a much more larger $j$. From the
reasoning described before, $f_p \to f$, and ${\cal F}(f)=f$,
$\langle f, 1\rangle=0$ and $f$ is differentiable. It is easy to
see that the modulus of concavity of $f$ is bounded (the iteration
by ${\cal F}$  does not decrease it). We can add a constant  to
$f$ and by linearity of ${\cal F}$ we also get a new fixed point
for ${\cal F}$ (note that ${\cal F}(1)=1$). Therefore, we can
assume $f= e^{-g}$ for some $g$.

In this way,  we obtain a contradiction with the uniqueness in
Theorem 26 in\cite{GV}.
\end{proof}

Suppose $\int g(x) \, \theta(x) dx =0$. For $\epsilon, h$ fixed,
then it follows from above that ${\cal F}^m (g) \to 0$ with
exponential velocity (according to the spectral gap).

Consider the backward stationary Markov  process $Y_n$
according to the transition $Q(x,v)$ and initial probability
$\theta$  as above.

\begin{theorem}
Given $f(x),g(x)$ with
$\int f(x) \, \theta(x) dx=\int g(x) \, \theta(x) dx =0$,
it follows
$$\int \, g(Y_0 ) \, f(Y_n)\, d \, P  \to 0,$$
with exponential velocity.
\end{theorem}
\begin{proof}

Note that
$$ \int g(Y_0 ) \, f(Y_1)\, d \, P= \int g(x) \, (\int Q(x,v) \, f(x-v) \, dv ) \,\theta(x)\,dx =$$
$$ \int\,  g(x)\,\,[\,  {\cal F}^*  (f)\, (x) \,]\,
\theta (x) dx \,=\, \int f(x) \, [\, {\cal F}(g) \,(x)\,] \,
\theta(x)\,  dx.
$$

In the same way,  for any $n$
$$\int \, g(Y_0 ) \, f(Y_n)\, d \, P=\, \int f(x) \, [\, {\cal F}^n (g) \,(x)\,] \, \theta(x)\,  dx .$$

The exponential decay of correlation follows from this.
\end{proof}

\begin{theorem}
Let $f(x),g(x)\in {\cal L}^2(\theta)$
be such that  $\int f(x) \, \theta(x) dx =\int g(x) \theta(x) dx=0$.
Then
$$\int \, g(X_0 ) \, f(X_n)\, d \, P  \to 0,$$
with exponential velocity.
\end{theorem}
\begin{proof}
Now, for analyzing the decay of  the forward system, $X_n$, with transition $\gamma(x,v)$, we have to consider the backwark operator
${\cal F}^*$, use the fact that  its exponential convergent, that is $({\cal F}^*)^n (g) \to 0$, if  $\int g(x) \, \theta(x) dx =0$,
and the result follows in the same way.
\vspace{1.0cm}
\end{proof}

\end{document}